\documentclass{amsart}
\usepackage[margin=1.0in]{geometry}
\usepackage{amsthm, amsmath, amssymb, mathtools, mathrsfs}
\usepackage{aliascnt}
\usepackage{hyperref}
\usepackage{breakurl}

\usepackage[textsize=footnotesize]{todonotes}

\usepackage{xcolor}
\definecolor{dblue}{rgb}{0,0,0.70}
\hypersetup{
	unicode=true,
	colorlinks=true,
	citecolor=dblue,
	linkcolor=dblue,
	anchorcolor=dblue
}

\makeatletter
\expandafter\g@addto@macro\csname th@plain\endcsname{%
	\thm@notefont{\bfseries}
}%
\expandafter\g@addto@macro\csname th@remark\endcsname{%
	\thm@headfont{\bfseries}
}%
\makeatother

\newtheorem{theorem}{Theorem}[section]	
\newtheorem*{theorem*}{Theorem}

\newaliascnt{lemma}{theorem}
\newtheorem{lemma}[lemma]{Lemma}
\aliascntresetthe{lemma}
\newtheorem*{lemma*}{Lemma}

\newaliascnt{proposition}{theorem}
\newtheorem{proposition}[proposition]{Proposition}
\aliascntresetthe{proposition}

\newaliascnt{corollary}{theorem}
\newtheorem{corollary}[corollary]{Corollary}
\aliascntresetthe{corollary}

\theoremstyle{remark}

\newaliascnt{remark}{theorem}

\aliascntresetthe{remark}
\newaliascnt{question}{theorem}
\newtheorem{question}[question]{Question}
\aliascntresetthe{question}

\newtheorem*{question*}{Question}

\newaliascnt{definition}{theorem}
\newtheorem{definition}[definition]{Definition}
\aliascntresetthe{definition}

\newaliascnt{example}{theorem}

\aliascntresetthe{example}

\usepackage[backend=biber, bibstyle=numeric, citestyle=numeric]{biblatex}
\renewbibmacro*{issue+date}{%
	\ifboolexpr{test {\iffieldundef{year}}
		and test {\iffieldundef{issue}}}
	{}
	{\printtext[parens]{%
			\printfield{issue}%
			\setunit*{\addspace}%
			\usebibmacro{date}}}%
	\newunit}

\NewBibliographyString{toappearin}
\NewBibliographyString{submittedto}
\DefineBibliographyStrings{english}{%
	toappearin  = {to appear in},
	submittedto = {submitted to},
}

\renewbibmacro*{in:}{%
	\ifboolexpr{not test {\iffieldundef{pubstate}}
		and (test {\iffieldequalstr{pubstate}{toappearin}}
		or test{\iffieldequalstr{pubstate}{submittedto}})}
	{\printtext{%
			\printfield{pubstate}\intitlepunct}%
		\clearfield{pubstate}}
	{\printtext{%
			\bibstring{in}\intitlepunct}}}

\newcommand{\lr}{\leftrightarrow}
\newcommand{\Ord}{\mathrm{Ord}}
\newcommand{\Fin}{\mathrm{Fin}}
\newcommand{\HF}{\mathrm{HF}}
\newcommand{\eps}{\mathrel{\mathtt{E}}}
\newcommand{\dom}{\operatorname{dom}}

\newcommand{\rank}{\operatorname{rank}}
\newcommand{\IE}{\mathsf{I}\mathcal{E}}
\newcommand{\SIE}{\mathsf{SI}\mathcal{E}}
\newcommand{\restricts}{\mathrel{\upharpoonright}}

\newcommand{\pair}{\mathsf{pair}}
\newcommand{\op}{\mathsf{op}}
\newcommand{\union}{\mathsf{union}}
\newcommand{\bininter}{\mathsf{bininter}}
\newcommand{\binunion}{\mathsf{binunion}}

\newcommand{\TC}{\operatorname{TC}}

\newcommand{\Dec}{\operatorname{Dec}}
\mathchardef\mhyphen="2D

\makeatletter
\def\oversortoftilde#1{\mathop{\vbox{\m@th\ialign{##\crcr\noalign{\kern3\p@}%
				\sortoftildefill\crcr\noalign{\kern3\p@\nointerlineskip}%
				$\hfil\displaystyle{#1}\hfil$\crcr}}}\limits}

\begin{document}
	
	\title{Constructive Ackermann's interpretation}
	\author{Hanul Jeon}
	\email{\href{mailto:hanuljeon95@gmail.com}{hanuljeon95@gmail.com}}
	\urladdr{\href{https://hanuljeon95.github.io}{https://hanuljeon95.github.io} }
	\dedicatory{This paper is dedicated to David C. McCarty.}
	\address{310 Malott Hall, Cornell University, Ithaca, New York 14853 USA}
	
	\subjclass[2010]{03F50, 03E70, 03F30, 03F65}
	\keywords{Heyting arithmetic, constructive set theory, finite set theory, interpretation}
	\thanks{This research was supported by BK21 SNU Mathematical Sciences Division.}
	
	\begin{abstract}
		The main goal of this paper is to formulate a constructive analogue of Ackermann's observation about finite set theory and arithmetic.
		We will see that Heyting arithmetic is bi-interpretable with $\mathsf{CZF^{fin}}$, the finitary version of $\mathsf{CZF}$. We also examine bi-interpretability between subtheories of finitary $\mathsf{CZF}$ and Heyting arithmetic based on the modification of Fleischmann's hierarchy of formulas, and the set of hereditarily finite sets over $\mathsf{CZF}$, which turns out to be a model of $\mathsf{CZF^{fin}}$ but not a model of finitary $\mathsf{IZF}$.
	\end{abstract}
	
	
	\maketitle
	
	\section{Introduction}
	Ackermann \cite{Ackermann1937} noticed in 1937 that $\mathsf{ZF}$ without Infinity is interpretable within $\mathsf{PA}$. Kaye and Wong \cite{KayeWong2007} improved Ackermann's result by showing that $\mathsf{ZF^{fin}}$, a finitary version of $\mathsf{ZF}$, is bi-interpretable with $\mathsf{PA}$.
	We may ask the same question for $\mathsf{HA}$: could we find a set theory that is bi-interpretable with $\mathsf{HA}$? One possible solution is to expand the language of set theory and add some axioms that are related to arithmetic, like \cite{Previale1994}. However, we want to stick our set-theoretic counterpart is closer to `standard' set theories as possible.
	
	Fortunately, we have a good start point: Aczel \cite{Aczel2013} characterized and gave a careful analysis of a weak subtheory $\mathsf{ACST}$ of constructive set theory $\mathsf{CZF}$, which is able to interpret Heyting arithmetic. Moreover, Aczel observed that there is a faithful interpretation from $\mathsf{HA}$ to $\mathsf{ACST}$. We may hope that we can find a set theory that is bi-interpretable with $\mathsf{HA}$ by extending $\mathsf{ACST}$, as Kaye and Wong capture set-theoretic counterpart of $\mathsf{PA}$ by extending $\mathsf{ZF}-\text{Infinity}$, which can be interpreted in $\mathsf{PA}$.
	
	It turns that our strategy works: we will see that $\mathsf{HA}$ is bi-interpretable with $\mathsf{ACST^{fin}}+\text{Set Induction}$, which is identical with $\mathsf{CZF^{fin}}$, the finitary version of $\mathsf{CZF}$. A bit surprisingly, Kaye and Wong's proof also works over constructive background, and we will see later how their proof can be implemented over a constructive setting.
	\begin{theorem*}\label{Theorem:MainResult}
		$\mathsf{HA}$ is bi-interpretable with $\mathsf{CZF^{fin}}$.
	\end{theorem*}
	Kaye and Wong also observed that their bi-interpretation is also a bi-interpretation between subtheories of $\mathsf{PA}$ and $\mathsf{ZF^{fin}}$, which is called $\mathsf{I\Sigma_n}$ and $\mathsf{\Sigma_n\mhyphen Sep}$ respectively, where $\mathsf{\Sigma_n\mhyphen Sep}$ is a theory of set theory that comprises Extensionality, Pairing, Empty set, Union, $\lnot$Infintiy, $\Delta_0$-Collection, $\Sigma_1\cup\Pi_1$-Set Induction and $\Sigma_n$-Separation.
	These subtheories rely on the hierarchy of formulas, which is given by prenex normal form, and there is little known about the complexity of predicate formulas over intuitionistic predicate logic, which renders finding appropriate constructive counterparts seemingly hard.
	However, we will see that modifying Fleischmann's hierarchy of formulas \cite{Fleischmann2010} yields a hierarchy on formulas of intuitionistic predicate logic, and also provides a nice subtheories of $\mathsf{HA}$ and $\mathsf{CZF^{fin}}$.
	As a result, we also have a constructive analogue of Kaye and Wong's analysis on subtheories:
	\begin{theorem*}
		For each $n\ge 1$, there are subtheories of $\IE_n$ of $\mathsf{HA}$ and $\SIE_n$ of $\mathsf{CZF^{fin}}$ that are bi-interpretable with each other.
		Moreover, adding the full law of excluded middle into $\IE_n$ and $\SIE_n$ yields $\mathsf{I\Sigma_n}$ and $\mathsf{\Sigma_n\mhyphen Sep}$ respectively.
	\end{theorem*}
	
	We know that $\mathsf{ZF}$ proves $\mathsf{ZF^{fin}}$ has a natural model, namely, the set of all hereditarily finite sets. We may ask whether a similar result holds $\mathsf{CZF}$. We will see the following theorem holds in \autoref{Section:Model}:
	\begin{theorem*}
		$\mathsf{CZF}$ proves the set of hereditarily finite sets $\HF$ exists and satisfies $\mathsf{CZF^{fin}}$. However, $\HF$ need not be a model of $\mathsf{IZF^{fin}}$.
	\end{theorem*}
	
	We develop facts that are necessary backgrounds in \autoref{Section:Interpretation} to \ref{Section:ordinalT}: we will define theories and interpretation between theories in \autoref{Section:Interpretation} and Heyting arithmetic in \autoref{Section:Heyting} respectively.
	We will define and discuss $\mathsf{ACST}$ and finitary constructive set theories in \autoref{Section:SetTheory}. We will also provide a streamlined version $\mathbb{T}$ of $\mathsf{CZF^{fin}}$ for further technical convenience.
	In \autoref{Section:Complexity}, we define subtheories of Heyting arithmetic and set theory based on a modification of Fleischmann's hierarchy of formulas \cite{Fleischmann2010}.
	\autoref{Section:interpretation} and \ref{Section:inverse} explain the interpretation \`a la Kaye and Wong is also a bi-interpretation between $\mathsf{HA}$ and $\mathsf{CZF^{fin}}$. It turns out in \autoref{Section:subtheory} that the interpretation can also interpret subtheories of $\mathsf{HA}$ and $\mathsf{CZF^{fin}}$ defined in \autoref{Section:Complexity}. In \autoref{Section:Model}, we will define and examine properties of the set of all hereditarily finite sets $\HF$ over $\mathsf{CZF}$.
	
	\section{Theories and interpretations}\label{Section:Interpretation}
	The central notion of this article is interpretation. Kaye and Wong \cite{KayeWong2007} and Aczel \cite{Aczel2013} adopted interpretations that are defined by Visser \cite{Visser2006}. In this article, we will follow Enderton's definition of interpretation (See Section 2.7 of \cite{Enderton2001},) with some adoptation of notions from \cite{KayeWong2007} and \cite{Aczel2013}.
	
	We will identify theories with a set of axioms. For given two theories $T_0$ and $T_1$, we write $T_0\vdash T_1$ if $T_0$ proves every axioms of $T_1$. If both $T_0\vdash T_1$ and $T_1\vdash T_0$ holds, we write $T_0\vdash\dashv T_1$ and say that $T_0$ and $T_1$ are \emph{identical}.
	
	To define an interpretation $\mathfrak{t}\colon T_0\to T_1$, we need formulas $\pi_\forall(x)$, $\pi_P(x_0,\cdots,x_{n-1})$ for $n$-ary predicate symbol $P$ and $\pi_f(x_0,\cdots,x_{n-1},y)$ for $n$-ary function symbol $f$.
	(We will regard constant symbols as nullary function symbols.)
	In addition, we assume that $T_1$ proves $\exists x \pi_\forall(x)$ and the functionality of $\pi_f$, that is,
	\begin{equation}
	T_1\vdash\forall x_0\cdots\forall x_{n-1}\exists ! y \pi_f(x_0,\cdots,x_{n-1},y).
	\end{equation}
	Then the interpretation $\mathfrak{t}$ sends a formula of $T_0$ to a formula of $T_1$ as follows:
	\begin{itemize}
		\item Let $s_0$, $\cdots$, $s_{n-1}$, $t_1$, $\cdots$, $t_m$ be terms, $f$ be a function symbol, and $P$ be a predicate symbol or $=$.
		Then $(P(f(s_0,\cdots,s_{n-1}),t_1,\cdots,t_m))^\mathfrak{t}$ is
		\begin{equation}
		\exists x_0\cdots\exists x_{n-1}\exists y
		\left[
		\bigwedge_{0\le i<n} (x_i=s_i)^\mathfrak{t}\land \pi_f(x_0,\cdots,x_{n-1},y)
		\land (P(y,t_1,\cdots,t_m))^\mathfrak{t}
		\right]
		\end{equation}
		We will apply the similar procedure if $f(t_i)$ appears on another argument of $P$. 
		\item $(P(x_0,\cdots,x_{n-1}))^\mathfrak{t}$ is $\pi_P(x_0,\cdots,x_{n-1})$, where each $x_i$ is a variable.
		\item $\mathfrak{t}$ respects logical connectives. For example, $(\phi\to\psi)^\mathfrak{t}$ is defined by $\phi^\mathfrak{t}\to\psi^\mathfrak{t}$.
		\item $(\forall x\phi(x))^\mathfrak{t}$ is $\forall x (\pi_\forall (x)\to \phi^\mathfrak{t}(x))$. $(\exists x\phi(x))^\mathfrak{t}$ is $\exists x (\pi_\forall (x)\land \phi^\mathfrak{t}(x))$. 
	\end{itemize}
	Every theory $T$ has the \emph{identity interpretation} $1_T$, which is defined by $\pi_\forall (x)\equiv (x=x)$, $\pi_P\equiv P$ and $\pi_f(\vec{x},y) \equiv (f(\vec{x})=y)$.
	For two interpretation $\mathfrak{s}, \mathfrak{t}\colon T_0\to T_1$ , we say they are the \emph{same} when $T_1\vdash \forall \vec{x} [\phi^\mathfrak{s}(\vec{x})\lr\phi^\mathfrak{t}(\vec{x})]$ for all formulas $\phi(\vec{x})$. In that case, we write $\mathfrak{s}=\mathfrak{t}$. The composition of two interpretation is the result of application of the two interpretations. 
	
	We can see that if $\mathfrak{s}\colon T_0\to T_1$ is an interpretation, then $T_0\vdash \phi$ implies $T_1\vdash \phi^\mathfrak{s}$. We call $\mathfrak{s}$ is \emph{faithful} if the converse also holds: that is, $\mathfrak{s}$ is faithful if $T_0\vdash \phi\iff T_1\vdash \phi^\mathfrak{s}$.
	For $\mathfrak{s}\colon T_0\to T_1$ and $\mathfrak{t}\colon T_1\to T_0$, $\mathfrak{s}$ and $\mathfrak{t}$ are the inverses of each other if $\mathfrak{ts}=1_{T_0}$ and $\mathfrak{st}=1_{T_1}$.	
	In that case, we call $\mathfrak{s}$ a \emph{bi-interpretation} between $T_0$ and $T_1$.
	
	\section{Heyting arithmetic}\label{Section:Heyting}
	Heyting arithmetic $\mathsf{HA}$ is the constructive counterpart of Peano arithmetic.
	There are various possible formulations of Heyting arithmetic: for example, we may take the language of arithmetic as the set of all primitive recursive functions and add axioms that define each primitive recursive functions.
	Since we want to analyze the relation between classical interpretation and constructive one, we choose the form given over the language $\mathcal{L}=\{0,S,+,\cdot\}$ with the following axioms:
	\begin{enumerate}
		\item $\forall x,y(Sx=Sy\to x=y)$,
		\item $\forall x (x=0\lor \exists y (x=Sy))$,
		\item $\forall x (x+0=x)$,
		\item $\forall x,y (x+Sy = S(x+y))$,
		\item $\forall x (x\cdot 0= 0)$,
		\item $\forall x,y (x\cdot Sy = x\cdot y + y)$,
		\item For each formula $\phi(x)$, $\phi(0)\land \forall x [\phi(x)\to\phi(S(x))]\to \forall x\phi(x)$.
	\end{enumerate}
	The last axiom is called the induction scheme. These set of axioms are strong enough to define primitive recursive functions and show they are provably total.
	Especially, we are interested in the totality of the exponential function, and we postulate it as an axiom $\mathsf{Exp}$. It is known that $\mathsf{HA}$ proves $\mathsf{Exp}$, but $\mathsf{Exp}$ could not be provable from a weaker subtheory of $\mathsf{HA}$.
	
	Heyting arithmetic does not include the law of excluded middle, but it proves the law of excluded middle for bounded formulas, that is, formulas whose quantifiers are of the form $\forall (x<y)$ or $\exists (x<y)$:
	\begin{proposition} \label{Proposition:BoundedLEM} \pushQED{\qed}
		Let $\phi$ be a bounded formula of $\mathsf{HA}$. Then $\mathsf{HA}$ proves $\phi\lor\lnot\phi$. \qedhere 
	\end{proposition}
	
	\section{Constructive finitary set theories}\label{Section:SetTheory}
	Aczel \cite{Aczel2013} defines an arithmetical version of constructive set theory $\mathsf{ACST}$ to analyze finite sets over constructive set theory $\mathsf{CZF}$. 
	We clarify some notions to define what $\mathsf{ACST}$ is. 
	A formula $\phi(x)$ of set theory is $\Delta_0$ if every quantifier in the formula is bounded, that is, every quantifier is of the form $\forall x(x\in a\to \cdots)$ or $\exists x(x\in a\land \cdots)$. We will abbreviate previous formulas into $\forall x\in a(\cdots)$ and $\exists x\in a(\cdots)$ respectively.
	
	Aczel defined $\mathsf{RCST}$ before defining $\mathsf{ACST}$.
	\begin{definition}
		$\mathsf{RCST}$ is the theory consisting of Extensionality, Empty set, Binary Intersection, Pairing and \emph{Global Union-Replacment Rule} $\mathsf{GURR}$: for each formula $\phi(u,v)$,
		\begin{equation}\label{Formula:GURR}
		\forall u \exists! v \phi(u,v)\to \forall x\exists y \forall z [z\in y\lr \exists v (z\in v \land \exists u\in x \phi(u,v))].
		\end{equation}		
	\end{definition}
	Intuitively, $\mathsf{GURR}$ states if $F$ is a class function and $x$ is a set, then $\bigcup F^"[x]$ is also a set. Hence $\mathsf{GURR}$ implies Union and Replacement.
	Some readers might wonder why $\mathsf{RCST}$ does not include a form of Separation. However, the following well-known result shows that $\mathsf{RCST}$ proves Separation for $\Delta_0$-formulas:
	
	\begin{lemma}\label{Proposition:BinaryIntersection} (Corollary 9.5.7 of \cite{AczelRathjen2010}) \pushQED{\qed}
		$\mathsf{BCST}$ without $\Delta_0$-Separation proves that $\Delta_0$-Separation is equivalent to Binary Intersection. \qedhere 
	\end{lemma}
	Here $\mathsf{BCST}$ comprises Extensionality, Pairing, Replacement, Union, Empty set, and $\Delta_0$-Separation. Since $\mathsf{RCST}$ proves every axiom of $\mathsf{BCST}$ except for $\Delta_0$-Separation, we can apply the previous lemma to $\mathsf{RCST}$.
	
	$\mathsf{ACST}$ is obtained by adding the induction scheme on natural numbers. We need to define what natural numbers are:
	\begin{definition}
		The class of natural numbers $\omega$ is defined as follows:
		\begin{equation}\label{Formula:DefinitonOmega}
		\omega = \{\alpha\in\Ord \mid \alpha^+ \subseteq \{0\}\cup \{\gamma^+\mid \gamma\in\Ord\}\},
		\end{equation}
		where $\alpha^+:=\alpha\cup\{\alpha\}$ and $\Ord$ is the class of all ordinals, that is, transitive sets whose elements are also transitive.
	\end{definition}
	
	\begin{definition}
		Mathematical Induction Axiom Scheme $\mathsf{MathInd}(\omega)$ is the following statement: for every definable class $X$, the following holds:
		\begin{equation}
		\operatorname{Ind}(X)\to \omega\subseteq X,
		\end{equation}
		where $\operatorname{Ind}(X)$ states $X$ is an inductive class:
		\begin{equation}
		\operatorname{Ind}(X)\equiv (0\in X \land \forall x\in X (x^+\in X) ).
		\end{equation}
		$\mathsf{ACST}$ is the theory obtained by adding $\mathsf{MathInd}(\omega)$ to $\mathsf{RCST}$.
	\end{definition}
	
	$\mathsf{ACST}$ is strong enough to do finitary mathematics. For example, the following theorem is a consequence of $\mathsf{ACST}$.
	\begin{lemma}[Primitive Recursion]\label{Lemma:PrimitiveRecursionOverOmega} ($\mathsf{ACST}$)
		Let $A$ and $B$ be classes and let $F_0\colon B\to A$, $F_1\colon B\times \omega\times A\to A$ be class functions. Then there is $H\colon B\times \omega\to A$ such that $H(b,0)=F_0(b)$ and $H(b,k+1) = F_1(b,k,H(b,k))$ for all $k\in\omega$.
	\end{lemma}
	\begin{proof}
		The proof of the lemma is available at Theorem 10.6 of \cite{Aczel2013}. We will give a direct proof that works over $\mathsf{ACST}$ for later analysis.
		
		A function $f$ is \emph{partially given under $b$ up to $m$} if $f\colon m^+\to A$ such that $f(0)=F_0(b)$ and $f(k^+)=F_1(b,k,f(k))$ for all $k\in m$.
		Take $\psi(b,k,f)$ if and only if $f$ is partially given under $b$ up to $m$ and $k\in\dom f$.
		
		We claim that $\forall b,m\exists f \psi(b,m,f)$ holds. 
		We will use induction on $m$. The case $m=0$ is obvious. Suppose that $\exists g \psi(b,m,g)$ holds. Take $u_0=F_1(b,k,g(m))$ and let $f=g\cup \{(m^+,u_0)\}$. Then $f$ witnesses $\psi(b,m^+,f)$.
		
		Now assume that $f_0$ and $f_1$ are partially given under $b$, up to $m_0$ and $m_1$ respectively.
		We can see that $k<\min(m_0,m_1)\to f_0(k)=f_1(k)$ by induction on $k$. Finally let
		\begin{equation}
		H(b,m)=x\iff \exists f [\text{$f$ is partially given under $b$ up to $m$ and }f(m)=x].\qedhere
		\end{equation}
	\end{proof}
	By Primitive recursion, we can see that $\mathsf{ACST}$ is able to define addition and multiplication over $\omega$.
	We will see later, however, that $\mathsf{ACST}$ does not suffice to be bi-interpretable with $\mathsf{HA}$. For example, $\mathsf{ACST}$ does not prove every set is finite, while the set theory simulated by $\mathsf{HA}$ seems should do, as $\mathsf{PA}$ does. 
	
	\begin{definition}
		Let $\Fin$ be the class of all finite sets, that is,
		\begin{equation}
		\Fin = \{x \mid \exists n\in\omega \exists f\colon n\to x (\text{$f$ is a bijection between $n$ and $x$})\}.
		\end{equation}
		The assertion $V=\Fin$ is that every set is finite, i.e., for each set $x$ there is $n\in \omega$ such that $x$ is a bijective image of $n$.
		One can show that $\mathsf{ACST}$ proves, for every finite set $x$, there is a unique such a natural number $n$. The readers can find its proof in Section 8.2 of \cite{AczelRathjen2010}. We call this natural number $n$ the \emph{cardinality} of $x$.
	\end{definition}
	Let $\mathsf{ACST^{fin}}$ be the theory $\mathsf{ACST}+(V=\Fin)$. It is easy to see that $\mathsf{ACST^{fin}}$ proves the axiom of choice
	\begin{equation}
	\forall a[\forall x\in a\exists y \phi(x,y)\to \exists f\in {^a}V \forall x\in a \phi(x,f(x))]
	\end{equation}
	and $\Delta_0$-excluded middle. Finally, let us consider the theory $\mathsf{ACST^{fin}}+\text{Set Induction}$.
	We will verify that the interpreted axiom is valid, but the verification takes effort per each axiom. Hence we prefer an axiom system as simple as possible. We will see that the following system $\mathbb{T}$ is a streamlined version of $\mathsf{ACST^{fin}}\mathrel{+}\text{Set Induction}$:
	\begin{definition}
		The theory $\mathbb{T}$  comprises the following axioms: Extensionality, Pairing, Union, Binary Intersection, Set Induction and $V=\Fin$.
	\end{definition}
	Obviously $\mathbb{T}$ is a subsystem of $\mathsf{ACST^{fin}}+\text{Set Induction}$. Moreover, we have
	\begin{proposition}\label{Proposition:AxiomofT}
		$\mathbb{T}$ and $\mathsf{ACST^{fin}}+\text{Set Induction}$ are identical.
	\end{proposition}
	\begin{proof}
		It suffices to show that $\mathsf{MathInd}(\omega)$ and $\mathsf{GURR}$ is derivable from $\mathbb{T}$.
		\begin{enumerate}
			\item $\mathsf{MathInd}(\omega)$:
			Let $X$ be an inductive class. Applying Set induction to the formula $x\in \omega\to x\in X$ yields the result. However, we will give an alternative proof for later analysis.
			
			Assume that $\alpha\in\omega$. Then we have $\alpha=0$ or $\alpha=\gamma^+$ for some ordinal $\gamma$. The former obviously implies $\alpha\in X$. 
			In the latter case, apply the set induction to $x\in \gamma^{++}\to x\in A$. Then we have $\alpha=\gamma^+\in A$.
			
			\item $\mathsf{GURR}$: Let $F$ be a class function. We will use the induction on the cardinality of $x$: assume that $\bigcup F^"[x]$ exists for sets $x$ of cardinality $x$. Then $\bigcup F^"[x\cup\{y\}] = \left(\bigcup F^"[x] \right)\cup F(y)$, whose existence follows from the listed axioms.\qedhere
		\end{enumerate}%
	\end{proof}
	
	Aczel showed that $\mathsf{ACST^{fin}}+\text{Set Induction}$ is identical with $\mathsf{CZF^{fin}}$, the finitary $\mathsf{CZF}$, which is obtained by replacing the axiom of Infinity in $\mathsf{CZF}$ to $V=\Fin$. The proof is direct, and it uses induction on the size of sets. We will give a proof of axiom of Strong collection from $\mathbb{T}$ for later analysis:
	\begin{proposition}\label{Propostion:StrColloverT}
		$\mathbb{T}$ proves the axiom of Strong Collection.
	\end{proposition}
	\begin{proof}
		Assume that $\forall x\in a\cup\{c\}\exists y \phi(x,y)$, where $a$ is a set of size $n$ and $c\notin a$. Assume inductively that Strong Collection holds for sets of size $n$, so we have $b$ such that
		$\forall x\in a\exists y\in b \phi(x,y)$ and $\forall y\in b\exists x\in a\phi(x,y)$.
		Take $d$ such that $\phi(c,d)$, then $b\cup \{d\}$ witnesses Strong Collection for $\phi$ and $a\cup\{c\}$.
	\end{proof}
	
	Kaye and Wong \cite{KayeWong2007} included the existence of transitive closure into their set-theoretic counterpart of $\mathsf{PA}$. Thus it is natural to ask whether our $\mathbb{T}$ should contain the existence of transitive closure as an axiom. The following lemma ensures it is provable from $\mathbb{T}$, so adding it is unnecessary:
	\begin{lemma}\label{Lemma:TransitiveClosure}
		$\mathbb{T}$ proves every set has a transitive closure; that is, for each $x$ there is a transitive set $y$ such that $x\subseteq y$ and for each transitive $z$ such that $x\subseteq z$, we have $y\subseteq z$. Furthermore, the class function $\TC(x)$ is definable. 
	\end{lemma}
	\begin{proof}
		We will show first that the \emph{class} $\TC(x)$ is uniformly definable. Consider $F_0(x)=x$ and $F_1(x,y)=x\cup\bigcup y$. By \autoref{Lemma:PrimitiveRecursionOverOmega}, there is $H$ such that $H(x,0)=x$ and $H(x,n+1)=x\cup\bigcup H(x,n)$. Now take $\TC(x)=\bigcup_{n\in\omega} H(x,n)$.
		
		We can show that $H(x,n)\subseteq y$ for all transitive $y$ such that $x\subseteq y$ by induction on $n$. Therefore $\TC(x)$ is the least transitive class that contains $x$. However, we do not know $\TC(x)$ is a set yet. We will show the following statement by induction on $x$:
		\begin{equation}\label{Formula:TCExistence}
		\exists u (\text{$u$ is transitive $\land$ $x\subseteq u\subseteq \TC(x))$}.
		\end{equation}
		Assume that \eqref{Formula:TCExistence} holds for all $y\in x$, i.e., for each $y\in x$ we can find a transitive $v$ such that $v\subseteq \TC(y)$.
		By Strong collection, we can find $c$ such that
		\begin{multline}
		[\forall y\in x\exists v\in c (\text{$v$ is transitive $\land$ $y\subseteq v\subseteq \TC(y))$}] \\ \land 
		[\forall v\in c\exists y\in x (\text{$v$ is transitive and $\land$ $y\subseteq v\subseteq \TC(y))$}].
		\end{multline}
		Now let $u=x\cup \bigcup c$. We claim that $u$ witnesses \eqref{Formula:TCExistence}. Since $u$ is a union of transitive sets, $u$ is transitive. $u\subseteq \TC(x)$ follows from $y\in x\to \TC(y)\subseteq \TC(x)$.
	\end{proof}
	
	We work over $\mathbb{T}$ in the remaining part of the section unless specified.
	We can see that the recursion theorem on sets holds since $\mathbb{T}$ proves every axiom of $\mathsf{CZF}$ except for Infinity.
	\begin{lemma}[Set Recursion]
		\label{Transfinite recursion}\label{Lemma:SetRecursion}
		Let $G$ be a total $(k+2)$-ary class function. Then there is a total $(k+1)$-ary class function $F$ such that
		\begin{equation}
		\forall \vec{x}\forall y [F(\vec{x},y) = G(\vec{x},y, \langle F(\vec{x},z)\mid z\in y\rangle)].
		\end{equation}
	\end{lemma}
	\begin{proof}
		We will follow the proof of Proposition 19.2.1 of \cite{AczelRathjen2010}. We present the whole proof for the sake of completeness and later analysis.
		
		Call $f$ be \emph{partially given under $\vec{x}$} if $f$ is a function of a transitive domain and $\forall y\in\dom f [f(y)=G(\vec{x},y,f\restricts y)]$.
		Take
		\begin{equation}
		\psi(\vec{x},y,f)\equiv (\text{$f$ is partially given under $\vec{x}$}) \land y\in\dom f.
		\end{equation}
		We will see that for given $\vec{x}$, $\forall y\exists f\psi(\vec{x},y,f)$ holds.
		We appeal to Set Induction on $y$: assume that $\forall u\in y\exists g \psi(\vec{x},u,g)$ holds. 
		By Strong Collection, we have a set $A$ such that $\forall u\in y\exists g\in A\psi(\vec{x},u,g)$ and $\forall g\in A\exists u\in y \psi(\vec{x},u,g)$.
		Let $f_0=\bigcup A$ and $u_0=G(\vec{x},y,\langle f_0(u)\mid u\in y\rangle)$. Take $f=f_0\cup\{(y,u_0)\}$.
		
		We want to claim that $f$ is a function. We need to ensure the following statement which can be shown easily by applying Set Induction: for any $g_0,g_1\in A$ and $x\in\dom(g_0)\cap\dom(g_1)$, we have $g_0(x)=g_1(x)$.
		Its upshot is that $f_0$ is a function. It is easy to see that the $\dom f_0$ is transitive, $y\subseteq\dom f_0$ and $\forall u\in\dom f_0[f(y)=G(\vec{x},y,f_0\restricts y)]$. Hence $\forall u\in \dom f_0 \ \psi(\vec{x},u,f_0)$.
		
		Now assume that $(a,b),(a,c)\in f$. Then ether both of them is a member of $f_0$ or one of them is a member of $\{(y,u_0)\}$.
		In the latter case, we can see that $b=c$ by applying $\psi(\vec{x},u,f_0)$ for all $u\in\dom f_0$ or the definition of $u_0$. Hence $f$ is a function. Checking the remaining conditions of $\psi(\vec{x},y,f)$ is direct, so we omit it.
	\end{proof}
	
	A class $\Phi$ is an inductive definition if it is a class of pairs. Every inductive class $\Phi$ is associated with a class of consequences $\Gamma_\Phi(A) = \{x\mid (A,x)\in\Phi\}$. A class $C$ is $\Gamma_\Phi$-closed if $\Gamma_\Phi(C)\subseteq C$.
	\begin{proposition}\label{Lemma:ClassInductiveDefinition} \pushQED{\qed}
		(Class Inductive Definition Theorem) Let $\Phi$ be an inductive definition.
		Then there is a least $\Gamma_\Phi$-closed class. \qedhere 
	\end{proposition}
	See Theorem 12.1.1 of \cite{AczelRathjen2010} for the proof of Class Inductive Definition Theorem. Note that the proof in \cite{AczelRathjen2010} works over $\mathsf{BCST}$ with Strong collection and Set induction, which are provable from $\mathbb{T}$.
	
	We will conclude this section with the following lemma, which asserts the composition of a $\Delta_0$-formula and a class function is still decidable.
	\begin{lemma} \label{Lemma:ClassFtnParDecidable}
		Let $F$ be a class function and $\phi(x)$ be a $\Delta_0$-formula. Then $\phi(F(x))\lor \lnot\phi(F(x))$ holds.
	\end{lemma}
	\begin{proof}
		It follows from the equivalence between $\phi(F(x))\lor \lnot\phi(F(x))$ and
		\begin{equation}
		\forall y [(y=F(x))\to (\phi(y)\lor\lnot\phi(y))].\qedhere
		\end{equation}
	\end{proof}
	
	\section{Complexity of formulas and subtheories}\label{Section:Complexity}
	We will analyze the interpretation over subtheories of $\mathbb{T}$ and $\mathsf{HA}$, as Kaye and Wong scrutinize subtheories of $\mathsf{PA}$ and finitary $\mathsf{ZFC}$.
	Kaye and Wong analyzes these theories along the complexity of Induction schemes and Separation. It calls into a need for a new hierarchy for arithmetical formulas, so we introduce a variation of Fleischmann's hierarchy of formulas \cite{Fleischmann2010}.
	\begin{definition}
		Let $\Phi$ and $\Psi$ be sets of formulas over the language of arithmetic or set theory.
		Then the set $\mathcal{E}(\Phi)$ is defined as the closure of $\Phi$ under $\land$, $\lor$, bounded quantifications and $\exists$.
		The set $\mathcal{U}(\Phi,\Psi)$ is defined as the smallest set such that
		\begin{enumerate}
			\item $\Phi\subseteq \mathcal{U}(\Phi,\Psi)$,
			\item $\mathcal{U}(\Phi,\Psi)$ is closed under $\land$, $\lor$, $\forall$, and bounded quantifications, and
			\item if $\psi\in\Psi$ and $\phi\in \mathcal{U}(\Phi,\Psi)$ then $\psi\to\phi$ is in $\mathcal{U}(\Phi,\Psi)$.
		\end{enumerate}
	\end{definition}
	\begin{definition}
		Let $\mathcal{E}_0=\mathcal{U}_0$ be the collection of all bounded formulas.
		For each $n\ge 1$, define $\mathcal{E}_n := \mathcal{E}(\mathcal{U}_{n-1})$ and $\mathcal{U}_n:=\mathcal{U}(\mathcal{E}_{n-1},\mathcal{E}_{n-1})$.
	\end{definition}
	Then $\mathcal{E}_n$ and $\mathcal{U}_n$ form increasing sequences of sets. 
	Moreover, we have the following fact by modifying the proof of Theorem 3.10 of \cite{Fleischmann2010}:
	\begin{proposition} \pushQED{\qed}
		$\bigcup_{n<\omega}\mathcal{E}_n = \bigcup_{n<\omega}\mathcal{U}_n$ and they are equal to the set of all formulas.
		Moreover, $\mathcal{E}_n$ and $\mathcal{U}_n$ are subsets of $\mathcal{E}_{n+1}$ and $\mathcal{U}_{n+1}$ respectively. \qedhere 
	\end{proposition}
	We will confuse $\mathcal{E}_n$ and $\mathcal{U}_n$ with a family of formulas that are \emph{provably} equivalent to an $\mathcal{E}_n$ and $\mathcal{U}_n$ formula respectively.
	Note that $\mathcal{E}_1$-formulas over the language of set theory are also known as $\Sigma$-formulas (cf. Definition 19.1.2 of \cite{AczelRathjen2010}.)
	$\Sigma$-reflection principle over $\mathsf{IKP}$ (Theorem 19.1.4 of \cite{AczelRathjen2010}) shows every $\mathcal{E}_1$-formula over $\mathsf{IKP}$ is in fact $\Sigma$-formula.
	Moreover, $\mathcal{E}_n$ and $\mathcal{U}_n$ classes are classically equivalent to $\Sigma_n$ and $\Pi_n$ classes respectively.
	
	We will analyze some definitions and theorems under the mentioned hierarchy. 
	\begin{definition}
		Let $\IE_n$ be a subtheory of $\mathsf{HA}$ where the full induction scheme is weakened to the induction scheme for $\mathcal{E}_n$-formulas. 
		$\SIE_n$ be a subtheory of $\mathbb{T}$ that restricts Set Induction schemes to $\mathcal{E}_n$-formulas.
	\end{definition}
	$\IE_1$ is still strong enough to show that every primitive recursive function is definable and total. Especially, $\IE_1$ proves $\mathsf{Exp}$.
	On the set-theoretic side, we can see that $\Delta_0\mhyphen\mathsf{LEM}$, the axiom of power set are still provable from $\SIE_1$.
	Moreover, we have the following results:
	\begin{lemma}[Mathematical Induction for $\mathcal{E}_n$-formulas] \label{Lemma:MathIndEn}
		Let $n\ge 1$. Then $\SIE_n$ proves $\operatorname{Ind}(X)\to \omega\subseteq X$ for each class $X$ that is given by an $\mathcal{E}_n$-formula.
	\end{lemma}
	\begin{proof}
		Note that we cannot apply Set Induction to $x\in\omega\to x\in X$ since the complexity of $x\in\omega$ is $\mathcal{E}_1$, so the complexity of the whole formula could not be $\mathcal{E}_n$. 
		Therefore, we instead apply the alternative proof of \autoref{Proposition:AxiomofT}: observe that $x\in \gamma^{++}\to x\in A$ is $\mathcal{E}_n$ if $A$ is an $\mathcal{E}_n$-class since it is equivalent to $\exists z\in \gamma^{++}(z=x\land z\in A)$, so the previous proof works.
	\end{proof}
	
	The following results can be obtained by modifying the proofs in \autoref{Section:SetTheory}:
	\begin{lemma} \pushQED{\qed}
		Let $n\ge 1$.
		Then the axiom of Strong Collection for $\mathcal{E}_n$-formulas are provable from $\SIE_n$. \qedhere
	\end{lemma}
	
	\begin{lemma}[Primitive recursion for $\mathcal{E}_n$-formulas]
		\label{Lemma:PrimitiveRecursionEn} \pushQED{\qed}
		$\mathsf{SI}\mathcal{E}_n$ proves the following:
		let $n\ge 1$ and $m\ge n$.
		Assume that $A$ and $B$ are $\mathcal{E}_m$-definable classes and take $\mathcal{E}_m$-definable class functions $F_0\colon B\to A$ and $F_1\colon B\times\omega\times A\to A$. Then there is a unique $H\colon B\times \omega\to A$ such that $H(b,0)=F_0(b)$ and $H(b,k+1)=F_1(b,k,H(b,k))$. \qedhere 	
	\end{lemma}
	Note that \autoref{Lemma:PrimitiveRecursionEn} shows the addition and multiplication on $\omega$ is still definable over $\SIE_n$.
	
	\begin{lemma}
		$\SIE_1$ proves the existence and $\mathcal{E}_1$-definability of $\TC(x)$.
	\end{lemma}
	\begin{proof}
		The proof of \autoref{Lemma:TransitiveClosure} still works over $\SIE_1$. 
		Moreover, $F_0$ and $F_1$ in the proof of \autoref{Lemma:TransitiveClosure} is $\mathcal{E}_1$, so $H$ and $\TC$ is also $\mathcal{E}_1$ by \autoref{Lemma:PrimitiveRecursionEn}.
	\end{proof}
	
	Hence we can carry on the usual proof of Set Recursion theorem and Class Inductive Definition Theorem for $\mathcal{E}_n$-formulas over $\mathsf{SI}\mathcal{E}_n$ for $n\ge 1$. The readers could consult with Section 9.3 and Chapter 12 of \cite{AczelRathjen2010} for the usual proof of Set Recursion Theorem and Class Inductive Definition Theorem.

	\begin{lemma}[Set Recursion for $\mathcal{E}_n$-formulas]
		\label{Lemma:SetRecursionEn} \pushQED{\qed}
		Assume $n\ge 1$. Then the following statement is provable from $\SIE_n$:
		Let $G$ be a total $(k+2)$-ary class function of complexity $\mathcal{E}_m$ for $m\ge n$.
		Then there is a total $(k+1)$-ary class function of complexity $\mathcal{E}_m$ such that $\forall \vec{x}\forall y [F(\vec{x},y)=G(\vec{x},y,\langle F(\vec{x},z)\mid z\in y \rangle)].$ \qedhere 
	\end{lemma}
	
	\begin{corollary}[Class Inductive Definition Theorem for $\mathcal{E}_n$-definitions]\label{Corollary:InductiveDefEn} \pushQED{\qed}
		Assume $n\ge 1$ and $m\le n$.
		Let $\Phi$ be an inductive definition of complexity $\mathcal{E}_m$. Then $\SIE_n$ proves there is a least $\Gamma_\Phi$-closed class, whose complexity is $\mathcal{E}_m$. \qedhere 
	\end{corollary}
	
	\section{Ordinals over $\SIE_1$}\label{Section:ordinalT}
	We will work over $\SIE_1$ in this section unless specified. The aim of this section is to prove the following theorem:
	\begin{theorem}\label{Theorem:ConsequenceofStrNegInf}
		The following statements hold over $\SIE_1$:
		\begin{enumerate}
			\item $\Ord=\omega$,
			\item There is a bijection from $V$ to $\Ord$.
		\end{enumerate}
	\end{theorem}
	
	We defined ordinals as transitive sets whose elements are also transitive, and $\Ord$ as the class of all ordinals. For example, every member of $\omega$ is an ordinal:
	\begin{lemma}
		Every member of $\omega$ is an ordinal.
	\end{lemma}
	\begin{proof}
		The proof can be done by induction on $n$.
	\end{proof}
	
	It is known that if every ordinal satisfies $\in$-least principle, then $\mathsf{\Delta_0\mhyphen LEM}$ holds. We will show the converse: that is, we will prove that every ordinal satisfies $\in$-least principle. In fact, we will see that some classical properties of ordinals hold. The proof is not difficult, but we present every detail of it. The reader should bear in mind that $\SIE_1$ proves the law of excluded middle for $\Delta_0$-formulas, so we may use it freely.
	
	\begin{lemma}
		If $\alpha$ is an ordinal and $A$ is an inhabited subset of $\alpha$, then $A$ has an $\in$-minimal element.
	\end{lemma}
	\begin{proof}
		Observe that the assertion `There is an $\in$-minimal element $a\in A$' is a $\Delta_0$ statement. Hence either $A$ has an $\in$-minimal element or every element of $A$ is not $\in$-minimal.
		Assume that the latter holds. We will prove $\forall a\in A (a\notin A)$ from the assumption by appealing to Set induction.
		Let $a$ be a set such that $b\notin A$ for all $b\in a$. If $a\in A$, then $a$ is an $\in$-minimal element of $A$, contradicting with the assumption on $A$. Therefore $a\notin A$. By Set induction, $a\notin A$ holds for all $a$. This contradicts with that $A$ is inhabited.
	\end{proof}
	
	\begin{lemma}
		If $\alpha$ and $\beta$ are an ordinals, then exactly one of $\alpha\in\beta$, $\alpha=\beta$ or $\alpha\ni\beta$ holds.
	\end{lemma}
	\begin{proof}
		We will use set induction on $\alpha$ and $\beta$ simultaneously.
		That is, assume that either $\gamma\in\delta$, $\gamma=\delta$ or $\gamma\ni\delta$ holds for all $\gamma\in\alpha$ and $\delta\in\beta$ if $\gamma$ and $\delta$ are ordinals.
		
		Assume that $\alpha$ and $\beta$ are ordinals. By $\Delta_0$-excluded middle, we have $\alpha=\beta$ or $\alpha\neq\beta$. If the latter holds, then
		\begin{equation}
		\lnot(\forall\gamma\in\alpha(\gamma\in\beta)\land \forall\gamma\in\beta(\gamma\in\alpha)),
		\end{equation}
		so either $\alpha\setminus\beta$ or $\beta\setminus\gamma$ has an element. Without loss of generality assume that $\gamma$ is an $\in$-minimal element of $\alpha\setminus\beta$.
		We want to show $\gamma=\beta$.
		
		Assume the contrary that $\gamma\neq\beta$ holds.
		Then one of $\gamma\setminus\beta$ or $\beta\setminus\gamma$ is inhabited. If $\delta\in\gamma\setminus\beta$, then $\delta\in\alpha\setminus\beta$ since $\delta\in\gamma\in\alpha$, contradicting with the minimality of $\gamma$.
		If $\delta$ is an $\in$-minimal element of $\beta\setminus\gamma$, then either $\gamma=\delta$, $\gamma\in\delta$ or $\gamma\ni\delta$ by inductive hypothesis. The two former hypotheses implies $\gamma\in\beta$, which contradicts with $\gamma\notin\beta$. The latter one contradicts with $\delta\in\beta\setminus\gamma$. In total, we have a contradiction.
		Therefore $\lnot(\gamma\neq\beta)$, so we have $\gamma=\beta$ due to $\Delta_0$-excluded middle.
	\end{proof}
	
	\begin{proof}[Proof of \autoref{Theorem:ConsequenceofStrNegInf}]
		We first claim that $\Ord=\omega$ holds. We will show by induction on $n$ that if there is an injection from $n$ to an ordinal $\alpha$, then $n\subseteq \alpha$.
		The case $n=0$ is trivial. Let assume the inductive hypothesis holds for $n$ and there is an injection $f\colon n+1\to \alpha$. Since we have an injection from $n$ to $\alpha$, $n\subseteq \alpha$. Moreover, there is $\beta\in\alpha$ which is different from members of $n$. Since both $n$ and $\beta$ are ordinals, we have $n\in\beta$, $n=\beta$ or $n\ni\beta$ holds. However, the latter case never happens, and the remaining two implies $n\in\alpha$.
		
		We will use $\Sigma$ and $\mathfrak{p}$ functions defined by \cite{KayeWong2007} to show $\Ord\cong V$.
		By \autoref{Lemma:PrimitiveRecursionOverOmega}, we can define $\hat{\Sigma}:\omega\times \mathcal{P}(\omega)\to\omega$ recursively as follows: $\hat{\Sigma}(0,x)=0$ and
		\begin{equation}
		\hat{\Sigma}(c+1,x) = 
		\begin{cases}
		\hat{\Sigma}(c,x), & \text{if }c+1\notin x,\\
		\hat{\Sigma}(c,x) + (c+1), & \text{if }c+1\in x,
		\end{cases}
		\end{equation}
		for all $c\in\omega$ and $x\in\mathcal{P}(\omega)$. Now take $\Sigma(x) = \hat{\Sigma}(\bigcup x, x)$. Finally, let
		\begin{equation}
		\mathfrak{p}(x) = \Sigma(\{2^{\mathfrak{p}(y)}\mid y\in x\}),
		\end{equation}
		where the exponentiation is a natural number exponentiation. We can see that $\hat{\Sigma}$ is $\mathcal{E}_1$, thus so does $\Sigma$ by \autoref{Lemma:PrimitiveRecursionEn}. \autoref{Lemma:SetRecursionEn} ensures 
		$\mathfrak{p}$ is well-defined and is $\mathcal{E}_1$. Moreover, we can show that $\mathfrak{p}:V\to\omega$ is a bijection as \cite{KayeWong2007} did: injectivity uses Set induction, and surjectivity uses induction on $\omega$.
	\end{proof}

	\section{The interpretation}\label{Section:interpretation}
	The main theorem of this paper is as follows:
	\begin{theorem}
		$\mathsf{HA}$ and $\mathbb{T}$ are bi-interpretable.
	\end{theorem}
	We first define an interpretation from $\mathbb{T}$ to $\mathsf{HA}$.
	Work over $\mathsf{HA}$, and define a primitive recursive binary relation $\eps$ as follows:
	\begin{equation}
	a\eps b \iff \exists r<2^a\exists m [b = (2m+1)\cdot 2^a + r]
	\end{equation}
	Intuitively, $a\eps b$ means the $a$th digit of the binary representation of $b$ is 1 as classically did. Since $m$ in the above formula satisfies $m<b$, $a\eps b$ is equivalent to a decidable formula.
	
	Define an interpretation $\mathfrak{a}\colon T\to \mathsf{HA}$ as follows: the domain of the interpretation is the whole natural numbers. Take $(a\in b)^\mathfrak{a} \equiv (a\eps b)$.
	We will show that every axiom of $\mathbb{T}$ is valid under $\mathfrak{a}$. The proof for $V=\Fin$ takes more effort, so we will postpone the proof of its validity.
	\begin{theorem}\label{Theorem:asendsTtoHA}
		If $\sigma$ is an axiom of $\mathbb{T}$ except for $V=\Fin$, then $\mathsf{HA}\vdash \sigma^\mathfrak{a}$.
	\end{theorem}
	\begin{proof}
		\begin{enumerate}
			\item Extensionality: The interpreted Extensionality states the following: if two natural numbers have the same binary representation, then they are the same. Since the existence and uniqueness of binary representation only relies on division algorithm and induction, which are still valid over $\mathsf{HA}$, the interpreted Extensionality is valid.
			
			\item Pairing: Consider the primitive recursive function $\pair(a,b)$ defined by $\pair(a,a) = 2^a$, and $\pair(a,b) = 2^a+2^b$ if $a\neq b$.
			It is easy to see that if $c=\pair(a,b)$, then $(\{a.b\} = c)^\mathfrak{a}$ holds.
			
			\item Union: Consider the following primitive recursive functions: define $\binunion(a,b)$ as $\binunion(a,0)=0$ and 
			\begin{equation}
			\binunion(a,2^c+b')  =
			\begin{cases}
			a & \text{if $b'=0$ and $c\eps a$},\\
			a + 2^c & \text{if $b'=0$ and $\lnot c\eps a$},\\
			\binunion(\binunion(a,2^c),\binunion(a,b')) & \text{otherwise.}
			\end{cases}
			\end{equation}
			for $b=2^c+b'$, $b'<2^c$. Now let $\union(0)=0$ and $\union(a) = \binunion(c,\union(a'))$ for $a=2^c+a'$, $a'<2^c$.
			Then we can see that $\union(\pair(a,b)) = \binunion(a,b)$ holds. Moreover, if $c=\union(a)$ then $(\bigcup a = c)^\mathfrak{a}$ holds.
			
			\item Binary Intersection: Define $\bininter(a,b)$ by primitive recursion on $b$ as follows: $\bininter(a,0)=0$ and
			\begin{equation}
			\bininter(a,2^c+b') = 
			\begin{cases}
			2^c & \text{if $b'=0$ and $c\eps a$},\\
			0 & \text{if $b'=0$ and $\lnot c\eps a$},\\
			\binunion(\bininter(a,2^c),\bininter(a,b')) & \text{otherwise.}
			\end{cases}
			\end{equation}
			for $b=2^c+b'$ and $b'<2^c$. Then $\bininter(a,b)$ witnesses the intersection of $a$ and $b$.
			
			\item Set induction: It directly follows from the induction of $\mathsf{HA}$ and the fact $a\eps b\to a<b$ for all $a$ and $b$.
		\end{enumerate}
	\end{proof}
	
	The case for $V=\Fin$ needs some preparation.
	We need an interpreted version of various notions to describe the interpreted axiom, so we define them. The ordered pair $\op(a,b)$ of $a$ and $b$ is $\pair(\pair(a,a),\pair(a,b))$.
	The von Neumann ordinal is defined recursively as follows: $v(0)=0$ and 
	\begin{equation}
	v(n+1)=\binunion(v(n),2^{v(n)}). 
	\end{equation}
	
	\begin{lemma}\label{Lemma:InterpretedvonNeumannCondition}
		If $(a\in\omega)^\mathfrak{a}$, then $a=v(n)$ for some $n$.
	\end{lemma}
	\begin{proof}
		We will use induction on $a$. If $a=0$, then take $n=0$.
		Assume that our theorem holds for all $c<a$, and $a$ satisfies $(a\in\omega)^\mathfrak{a}$. Since $a>0$, there is $\gamma$ such that
		\begin{equation}\label{Formula:SuccessorAssumption}
		\exists c [(c \in\Ord)^\mathfrak{a} \land a = \binunion(c,2^c)]
		\end{equation}
		We can see that $a= \binunion(c,2^c)$ implies $c<a$. By the inductive hypothesis, $c=v(n)$ for some $n$. By definition of $v$, we have $a=v(n+1)$.
	\end{proof}
	
	\begin{theorem}
		$\mathsf{HA}$ proves $(V=\Fin)^\mathfrak{a}$.
	\end{theorem}
	\begin{proof}
		Note that the word `function' in this proof means a binary relation with the definining condition of a function. Before to describe the proof, let me define a \emph{size} $\sigma(a)$ of a natural number: $\sigma(0)=0$ and $\sigma(a) = 1+\sigma(a')$ if $a=2^c+a'$ for some $c<a$ and, $a'<2^c$.
		
		We will use induction on $a$. If $a=0$, then 0 witnesses the bijection between $0$ and $v(0)=0$.
		Let $a=2^c+a'$ for $a'<2^c$. Assume inductively that we have a function $f'$ from $a'$ to $v(\sigma(a'))$, and the inverse $g'$ of $f'$. We claim that the relation
		\begin{equation}
		f = \binunion(f',2^{(\op(c, v(\sigma(a')))})
		\end{equation}
		is a function from $a$ and $v(\sigma(a)) = v(\sigma(a')+1)$ with the inverse function
		$g = \binunion(g',2^{\op(v(\sigma(a')),c)})$
		
		Since $f'$ is a function of domain $a'$ and $a'<2^c$, the domain of $f'$ does not contain $c$. Hence $f$ is a function. It is obvious that the domain of $g'$, namely $v(\sigma(a'))$, does not contain $v(\sigma(a'))$. It shows $g$ is a function.
		It remains to show that $f$ and $g$ are inverses of each other, but it is clear from the properties of $f'$ and $g'$ and the definition of $f$ and $g$.
	\end{proof}
	
	In summary, we have
	\begin{corollary}\pushQED{\qed}
		$\mathfrak{a}$ is an interpretation from $\mathbb{T}$ to $\mathsf{HA}$.\qedhere 
	\end{corollary}
	
	\section{The inverse interpretation}\label{Section:inverse}
	We follow the inverse interpretation given by \cite{KayeWong2007}: we first define the ordinal interpretation $\mathfrak{o}$ from $\mathsf{HA}$ to $\mathbb{T}$, and compose it with $\mathfrak{p}$.
	
	\begin{definition}
		The ordinal interpretation $\mathfrak{o}$ is defined as follows: $0^\mathfrak{o}$ is the empty set, $S^\mathfrak{o}(x) = x\cup\{x\}$. Interpretation of addition and multiplication employs the corresponding operation on ordinals.
	\end{definition}
	Then we have
	\begin{theorem}
		$\mathfrak{o}$ is an interpretation from $\mathsf{HA}$ to $\mathbb{T}$.
	\end{theorem}
	\begin{proof}
		The only remaining axiom we need to check its validity is the induction scheme, and it follows from induction on $\omega$.
	\end{proof}
	
	$\mathfrak{o}$ is not an inverse interpretation of $\mathfrak{a}$ because of the interpretation of quantifiers. The formula $\forall x \phi(x)$ over $\mathsf{HA}$ is interpreted into $\forall x\in \omega \phi^\mathfrak{o}(x)$ under $\mathfrak{o}$, and its interpretation under $\mathfrak{a}$ is $\forall x [(x\in\omega)^\mathfrak{a}\to \phi^{\mathfrak{ao}}(x)]$, which is not equivalent to the original formula even if $\phi(x)$ is atomic.
	Kaye and Wong resolved this problem by relying on the bijection $\mathfrak{p}\colon V\cong \omega$.
	
	\begin{definition}
		The interpretation $\mathfrak{b}$ is defined as follows: the domain of the interpretation is $x=x$. If $t(\vec{x})$ is a term of $\mathsf{HA}$ whose variables are all expressed, then $t^\mathfrak{b}$ is defined as $t^\mathfrak{o}(\mathfrak{p}(\vec{x}))$.
	\end{definition}
	Then $\mathfrak{b}$ is an interpretation from $\mathsf{HA}$ to $\mathbb{T}$.
	Moreover, $\mathfrak{a}$ and $\mathfrak{b}$ are inverses of each other:
	\begin{theorem}\label{Theorem:BiInterpretation}
		$\mathfrak{ab}=1_{\mathsf{HA}}$ and $\mathfrak{ba}=1_\mathbb{T}$.
	\end{theorem}
	\begin{proof}
		Both $\mathfrak{ab}$ and $\mathfrak{ba}$ respects equality and domain, so it suffices to see that both interpretations preserve atomic symbols.
		
		To prove $\mathfrak{ba}$ is the identity, it is sufficient to see that $(x\in y)^\mathfrak{ba}$ is equivalent to $x\in y$. Note that $(x\in y)^\mathfrak{ba}$ is equivalent to $(x\eps y)^\mathfrak{b}$, in other words,
		\begin{equation}\label{Formula:inba}
		\exists r<2^{\mathfrak{p}(x)}\exists m\in\omega [ \mathfrak{p}(y) = (2m+1)2^{\mathfrak{p}(x)} + r]. 
		\end{equation}
		The above formulation $(x\in y)^\mathfrak{ba}$ is not a $\Delta_0$-formula for the following reasons: First, we do not know $\omega$ is a set, and it could be a proper class. Second, the definition of $\mathfrak{p}$ is not $\Delta_0$. Despite that, $(x\in y)^\mathfrak{ba}$ is a decidable formula by the fact that $x\eps y$ is equivalent to a decidable formula and \autoref{Lemma:ClassFtnParDecidable}.
		
		Let $y=\{z_k \mid k<n\}$ be an enumeration of $y$ such that $\mathfrak{p}(z_0)>\cdots>\mathfrak{p}(z_{n-1})$.
		We can see that for each $k<n$, there are $m$ and $r<\mathfrak{p}(z_k)$ such that $\mathfrak{p}(y) = (2m+1)\mathfrak{p}(z_k) + r$ by induction on $k$ and Euclidean division algorithm.
		If $x\in y$, then $x=z_k$ for some $k$, thus it satisfies \eqref{Formula:inba}.
		If $x\notin y$, then $\mathfrak{p}(x)$ is equal to none of the $\mathfrak{p}(z_k)$.
		By dividing cases, we can see that $\mathfrak{p}(y) = 2m 2^{\mathfrak{p}(x)} + r$ for some $m$ and $r<2^{\mathfrak{p}(x)}$. By uniqueness of the remainder and quotient, we have the negation of \eqref{Formula:inba}.
		Since both $x\in y$ and \eqref{Formula:inba} are decidable, we have the equivalence of these two formulas.
		
		It remains to show that $\mathfrak{ab}$ is the identity. It requires a sequel of lemmas on interpreted notions of $\mathbb{T}$.
		We can see that $(a\in\omega)^\mathfrak{a}$ \emph{if and only if} $a=v(n)$ for some $n$ by \autoref{Lemma:InterpretedvonNeumannCondition} and an easy inductive argument.
		Moreover, we can show the following fact by induction on $y$:
		\begin{lemma}\pushQED{\qed}
			For each $x$ and $y$, we have $(S^\mathfrak{a}(v(y))=v(y+1)$, 
			$v(x)+^\mathfrak{a}v(y)=v(x+y)$, $v(x)\cdot^\mathfrak{a} v(y)=v(x\cdot y)$, and
			$(v(x)^{v(y)})^\mathfrak{a} = v(x^y)$.\qedhere
		\end{lemma}
		Here the functions under $\mathfrak{a}$ are set-theoretic functions for von Neumann ordinals, and the functions appearing in the argument of $v$ are functions of the language of $\mathsf{HA}$. From this lemma, we have
		\begin{lemma}
			$\mathfrak{p}^\mathfrak{a}(x)=v(x)$; that is, $(\mathfrak{p}(x)=y)^\mathfrak{a}$ if and only if $v(x)=y$.
		\end{lemma}
		The proof uses induction on $x$: assume that the desired equality holds for all $y<x$.
		Let $x=2^c+x'$ for $x<2^c$. Then $2^c+x' = \binunion(2^c,x')$, so $(x=\{c\}\cup x')^\mathfrak{a}$. Hence
		$(\mathfrak{p}(x) = \mathfrak{p}(x'\cup\{c\}) = 2^{v(c)} + v(x'))^\mathfrak{a}$, and the previous lemma ensures $\mathfrak{p}^\mathfrak{a}(x) = v(x)$.
		
		It remains to check that function symbols $S$, $+$ and $\cdot$ are preserved under $\mathfrak{ab}$. We will only see the proof for $S$, as the remaining cases are analogous: $(S(x)=y)^\mathfrak{ab}$ if and only if $(\mathfrak{p}(x)\cup\{\mathfrak{p}(x)\}=\mathfrak{p}(y))^\mathfrak{a}$, which turns out to be equivalent to $\binunion(v(x),2^{v(x)}) = v(y)$. Hence $v(S(x))=v(y)$, which is equivalent to $S(x)=y$.
	\end{proof}
	
	We will conclude this section by correcting a result of Aczel.
	Aczel stated in Section 11 of \cite{Aczel2013} that the ordinal interpretation $\mathfrak{o}\colon \mathsf{HA}\to T$ is a faithful interpretation if $T$ is a subtheory of $\mathsf{IZF^{fin}}$ that contains $\mathsf{ACST}$.
	We know that $\mathsf{IZF^{fin}}$ is a classical theory, so it proves the interpreted version of semi-classical principles like $\mathsf{WLEM}$. However, we know that $\mathsf{HA}$ does not prove $\mathsf{WLEM}$. Hence there is no faithful interpretation from $\mathsf{HA}$ to $\mathsf{IZF^{fin}}$. The result holds, fortunately, if we correct $\mathsf{IZF^{fin}}$ to $\mathsf{CZF^{fin}}$:
	\begin{proposition}
		Let $T$ be a subtheory of $\mathsf{CZF^{fin}}$ that contains $\mathsf{ACST}$ as a subtheory. Then $\mathfrak{o}\colon \mathsf{HA}\to T$ is a faithful interpretation.
	\end{proposition}
	\begin{proof}
		We can see that $\mathsf{ACST}$ is capable of defining $\mathfrak{o}$, hence $\mathfrak{o}\colon \mathsf{HA}\to T$ is an interpretation.
		We will see that the composition $\mathfrak{ao}\colon  \mathsf{HA}\to\mathsf{HA}$ is faithful:
		Assume that $\mathsf{HA}\vdash \phi^\mathfrak{ao}$.
		By applying inverse interpretation $\mathfrak{b}$, we have $\mathsf{CZF^{fin}}\vdash \phi^\mathfrak{o}$. We can see that every quantifier of $\phi^\mathfrak{o}$ is of the form $\forall x\in\omega$ or $\exists x\in\omega$.
		We can see that the formula is still provable if we replace every variable $x$ to $\mathfrak{p}(x)$ and omit $\in\omega$ from quantifiers since $\mathfrak{p}$ is a definable bijection between $V$ to $\omega$. Hence $\mathsf{CZF^{fin}}\vdash \phi^\mathfrak{b}$, and we have $\mathsf{HA}\vdash \phi$. The main result follows directly from the previous argument.
	\end{proof}
	
	\section{Interpretating subtheories of $\mathsf{HA}$}\label{Section:subtheory}
	Kaye and Wong gave not only a bi-interpretation between $\mathsf{PA}$ and its set-theoretic counterpart, but also a bi-interpretation between their subtheories. It can be given by asserting the previous proof works over a subtheory of $\mathsf{PA}$ and finitary set theory.
	We can do a similar work over $\mathsf{HA}$ and $\mathbb{T}$ under the $\mathcal{E}_n$-hierarchy of formulas.
	
	Let us analyze the definitional complexity of $\mathfrak{a}$. As \cite{KayeWong2007} pointed out, the definition of $\mathfrak{a}$ only involves bounded formulas except for exponentiation function. However, the exponentiation function is definable by a $\mathcal{E}_1$-formula with the help of $\mathcal{E}_1$-induction. In fact, we can show the following general theorem:
	\begin{proposition}\label{Proposition:E1RecursionArithmetic} \pushQED{\qed}
		Every primitive recursive function is definable by an $\mathcal{E}_1$-formula over $\IE_1$. \qedhere
	\end{proposition}
	The proof of \autoref{Proposition:E1RecursionArithmetic} is the same with the classical counterpart of the theorem: namely, every primitive recursive function is $\Sigma_1$-definable over $\mathsf{I\Sigma_1}$.
	As a speical case of \autoref{Proposition:E1RecursionArithmetic}, we can see $a\eps b$ is definable by a $\mathcal{E}_1$-formula.
	Hence $\mathfrak{a}$ sends each $\mathcal{E}_n$ (or $\mathcal{U}_n$) formulas of set theory to $\mathcal{E}_n$ (or $\mathcal{U}_n$) formulas of arithmetic. From this, we can infer the following theorem:
	
	\begin{theorem}\pushQED{\qed}
		Let $n\ge 1$. Then $\mathfrak{a}$ is an interpretation from $\SIE_n$ to $\IE_n$. \qedhere
	\end{theorem}
	
	The case of $\mathfrak{b}$ involves with the complexity of $\mathfrak{o}$ and $\mathfrak{p}$.
	We can see that the definitional complexity of $\hat{\Sigma}$, $\mathfrak{p}$ and primitive recursive functions are $\mathcal{E}_1$ by \autoref{Lemma:PrimitiveRecursionEn} and \autoref{Lemma:SetRecursionEn}.
	As the definition of $\mathfrak{b}$ employs $\mathfrak{p}$, $\mathfrak{b}$ could increase the complexity of formulas.
	Fortunately, we can see that $\mathcal{E}_n$-classes are stable under this substitution since $\phi(F(x),y)$ is equivalent to $\exists z(z=F(x)\land \phi(z,y))$:
	\begin{lemma} ($\mathsf{BCST}$) \label{Lemma:EnPreservationbyEnFtn}\pushQED{\qed}
		Let $\phi(x,\vec{y})$ be an $\mathcal{E}_n$-formula for $n\ge 1$ and $F(x)$ be a $\mathcal{E}_1$-class function. Then $\phi(F(x),\vec{y})$ is an $\mathcal{E}_n$-formula. \qedhere 
	\end{lemma}
	
	Hence we have the following theorem:
	\begin{theorem}\label{Theorem:bInSubtheory}\pushQED{\qed}
		$\mathfrak{b}$ is an interpretation from $\IE_n$ to $\SIE_n$.\qedhere
	\end{theorem}
	
	We can see that the proof of \autoref{Theorem:BiInterpretation} works over $\IE_n$ and $\SIE_n$. Therefore, we have the following:
	\begin{corollary}\pushQED{\qed}
		Let $n\ge1$. Then the interpretations $\mathfrak{a}\colon \SIE_n\to \IE_n$ and $\mathfrak{b}\colon \IE_n\to\SIE_n$ are inverses of each others.\qedhere
	\end{corollary}
	
	We will finish this section by showing that the previous results are exactly the constructive counterpart of \cite{KayeWong2007}.
	Kaye and Wong proved that $\mathfrak{a}$ and $\mathfrak{b}$ is bi-interpretations between $\mathsf{I\Sigma_n}$ and $\mathsf{\Sigma_n\mhyphen Sep}$.
	It is easy to see that $\IE_n$ with the full law of excluded middle is $\mathsf{I\Sigma}_n$. The following theorem shows $\SIE_n$ is the constructive counterpart of $\mathsf{\Sigma_n\mhyphen Sep}$:
	\begin{theorem}
		Let $n\ge1$. Then $\SIE_n$ with the full law of excluded middle is identical with $\mathsf{\Sigma_n\mhyphen Sep}$.
	\end{theorem}
	\begin{proof}
		Let $\SIE_n^c$ be $\SIE_n$ with the full law of excluded middle. Since we have the full law of excluded middle, $\mathcal{E}_n$ and $\mathcal{U}_n$-classes coincide with $\Sigma_n$ and $\Pi_n$-classes respectively. We will divide the proof into two parts:
		
		\begin{lemma}
			$\SIE_n^c\vdash\mathsf{\Sigma_n\mhyphen Sep}$.
		\end{lemma}
		\begin{proof}
			Note that $\SIE_n^c\vdash \Sigma_n$-Strong Collection by the proof of \autoref{Propostion:StrColloverT}.
			We can see that $\Sigma_n$-Strong Collection with the full law of excluded middle proves $\Sigma_n$-Separation. Proving $\Delta_0$-Collection and $\lnot$Infinity from $\SIE_n^c$ are easy. In sum, we have shown that $\SIE_n^c$ proves every axioms of $\mathsf{I\Sigma_n}$ except for $\Sigma_1\cup\Pi_1$-Induction.
			
			Proving the induction scheme is trivial if $n\ge 2$, since every $\Pi_1$-formula is $\Sigma_2$. However, there is another way to prove $\Pi_1$-induction from the remaining axioms, which works even if we only have $\mathcal{E}_1$-induction:
			Let $\phi(x)$ be a $\Pi_1$-formula. Assume the contrary that $\phi(x)$ does not satisfy set induction, in other words, $(\forall y\in x \phi(y))\to \phi(x)$ holds for all $x$, but $\lnot\phi(a)$ also holds for some $a$.
			Let $n$ be a cardinality of $\TC(a)$, and consider the set
			\begin{equation}
			X=\{m\in n^+ \mid \exists x [\lnot\phi(x) \land \exists f\colon m\to \TC(x) \text{ such that $f$ is a bijection.}]\}
			\end{equation}
			($\Sigma_1$-Separation is needed to define $X$.) Then $X$ is an inhabited set of ordinals, so it has a minimal element $m$. 
			Let $x$ be a set such that $\lnot\phi(x)$ and the cardinality of $\TC(x)$ is $m$.
			For each $y\in\TC(x)$, we have $|\TC(y)|<|\TC(x)|$, so $\phi(y)$. Therefore $\phi(y)$ for all $y\in x$, and it implies $\phi(x)$, a contradiction.
		\end{proof}
		
		\begin{lemma}
			$\mathsf{\Sigma_n\mhyphen Sep}\vdash\SIE_n^c$.
		\end{lemma}
		\begin{proof}
			Observe that $\Delta_0$-Set Induction proves the axiom of Regularity.
			We will see that $\Sigma_n$-induction scheme holds over $\mathsf{I\Sigma_n}$.
			
			Assume that $\phi(x)$ is a $\Sigma_n$-formula that violates Set Induction, so that $\forall y\in x\phi(x)\to\phi(x)$ for all $x$, but there is $a$ such that $\lnot\phi(a)$.
			Now consider the set $X=\{x\in a \mid \lnot\phi(x)\}$. If $X$ is empty, then $\phi(a)$. If $X$ is not empty, choose an $\in$-minimal element $x$ of $X$, which would satisfy $\phi(x)$. In either case, we have a contradiction. Therefore $\phi(x)$ follows Set Induction scheme.
			
			Finally, the negation of Infinity implies $V=\Fin$ classically: 
			Assume that $V\neq \Fin$ holds. Let $a$ be a set that is not bijectable with any $n\in\omega$.
			Since we have $\Sigma_1$-Set Induction, we can apply \autoref{Lemma:SetRecursionEn} to define rank function as follows:
			\begin{equation}
			\rank x = \bigcup\{\rank y + 1 \mid y\in x\}.
			\end{equation}
			Note that we proved \autoref{Lemma:SetRecursionEn} over $\SIE_1$, so there is a possibility that we are using $V=\Fin$ in the proof of \autoref{Lemma:SetRecursionEn}. However, we do not need $V=\Fin$ in this proof of \autoref{Lemma:TransitiveClosure} and \autoref{Lemma:SetRecursionEn} if we have $\Sigma_1$-Collection, which is provable from $\mathsf{\Sigma_1\mhyphen Sep}$.
			
			We can show that if every element of $x$ has a rank smaller than $n$, then $x$ is finite: it follows from that $V_n$ exists for each natural number $n$, which is a theorem of $\mathsf{\Sigma_1\mhyphen Sep}$. Hence the set $X=\{\rank y\mid y\in a\}$ is infinite. Since $\omega\subseteq \bigcup X$, $\omega$ is a set. Therefore, we have the axiom of Infinity.
		\end{proof}
		
		Combining these two lemmas, we have $\mathsf{\Sigma_n\mhyphen Sep}\vdash\dashv\SIE_n^c$
	\end{proof}
	
	Note that the only properties of $\mathcal{E}_n$ used in this section for proving $\mathfrak{a}$ and $\mathfrak{b}$ are well-defined and bi-interpretations of each other are that $\mathcal{E}_n$ contains bounded formulas and $\mathcal{E}_n=\mathcal{E}(\mathcal{E}_n)$.
	Thus we can extend our argument to any class of formulas with certain conditions. We state it without proof:
	\begin{proposition}\pushQED{\qed}
		Let $\Gamma$ and $\Gamma'$ be collection of formulas over set theory and arithmetic respectively, such that both of $\Gamma$ and $\Gamma'$ contains $\Delta_0$.
		Assume that $\Gamma$ and $\Gamma'$ satisfies $\mathcal{E}(\Gamma)=\Gamma$ and $\mathcal{E}(\Gamma')=\Gamma'$. Furthermore, assume that $\mathfrak{a}$ sends $\Gamma$-formulas to $\Gamma'$-formulas and \textit{vice versa} for $\mathfrak{b}$.
		
		Then $\mathfrak{a}$ is a bi-interpretation between $\mathsf{SI}\Gamma$, a theory obtained by restricting Set Induction in $\mathbb{T}$ to $\Gamma$-formulas, and $\mathsf{I}\Gamma'$, which is obtained by restricting induction scheme of $\mathsf{HA}$ to $\Gamma'$-formulas. \qedhere 
	\end{proposition}
	
	\section{A natural model of $\mathsf{CZF^{fin}}$: the set of hereditarily finite sets}\label{Section:Model}
	In this section, we will work over $\mathsf{CZF}$ unless stated otherwise.
	
	$\mathsf{ZF}$ proves that the set of sets of all finite rank $V_\omega$ is a model of $\mathsf{ZF^{fin}}$. Moreover, we may regard $V_\omega$ as a natural model of $\mathsf{ZF^{fin}}$ since it is countable and every set given by Ackermann's intepretation falls into $V_\omega$.
	We may ask $\mathsf{CZF}$ can prove the existence of a natural model of $\mathsf{CZF^{fin}}$.
	One candidate is $V_\omega$, but we will not consider it for the following reasons: first, $V_\omega$ could not be a set in $\mathsf{CZF}$. (In fact, even the power set of 1 need not be a set in $\mathsf{CZF}$.) Second, we cannot ensure $V_\omega$ need not be countable even if we assume the axiom of power set. In fact, McCarty showed that $V[\mathcal{K}_0]$, the model given by Kleene realizability satisfies $\mathsf{IZF}$ with `$\mathcal{P}(1)$ is not subcountable.' (See Corollary 3.8.3 of \cite{McCartyPhD}.)
	Aczel \cite{Aczel2013} also introduced the class of hereditarily finite sets $\HF$, and provide some properties of it. We have to add some remarks on finiteness and related concepts over $\mathsf{CZF}$ before we can say about what $\HF$ is.
	
	Finite sets are sets that have a bijection with a von Neumann natural number. However, this notion of finite sets over constructive set theory is not well-behaved unlike classical finite sets. For example, a subset of a finite set need not be finite.
	It does not mean we have to take another definition of finite sets. Instead, we divide notions of finiteness:
	\begin{definition}
		A set $x$ is \emph{finite} if there is a bijection from $n\in\omega$ to $x$.
		$x$ is \emph{finitely enumerable}  if there is a surjection from $n\in\omega$ to $x$. 
		$x$ is \emph{subfinite} if $x$ is a subset of a finite set.
	\end{definition}
	
	The following condition has a essential role to characterize which subsets of a finite is again finite:
	\begin{definition}
		Let $A$ be a set. A subset $B\subseteq A$ is \emph{decidable} if $x\in B\lor x\notin B$ for all $x\in A$. $A$ is \emph{discrete} if the equality relation is a decidable subset of $A\times A$.
	\end{definition}
	
	\begin{proposition} (Proposition 8.1.11 of \cite{AczelRathjen2010}, $\mathsf{CZF}$) \pushQED{\qed}
		A set $x$ is finite if and only if $x$ is finitely enumerable and discrete. \qedhere
	\end{proposition}
	
	\begin{lemma}\label{Lemma:DecidableFinite}
		Every decidable subset of a finite set is finite.
	\end{lemma}
	\begin{proof}
		It suffices to show that every decidable subset of $n\in\omega$ is finite.
		We can show by induction on $n$ that if $\phi(m)$ is decidable on $n$ then either $\exists m\in n \phi(m)$ or $\forall m\in n \lnot\phi(m)$.
		Hence if $x\subseteq n$ is decidable, then $x$ is empty or inhabited.
		If $x$ is inhabited and $m_0\in x$, then the function
		\begin{equation}
		f(m)=\begin{cases}
		m&\text{if }m\in x,\\
		m_0&\text{otherwise}
		\end{cases}
		\end{equation}
		enumerates elements of $x$. Hence $x$ is finitely enumerable. Since $x\subseteq n$ is discrete, $x$ is finite.
	\end{proof}
	
	\begin{definition}
		Let $\Phi$ be an inductive definition given as $\Phi_\text{fin}:=\{(\{a\},a) \mid a\text{ is finite}\}$. Then $\HF$ is the least $\Gamma_{\Phi_\text{fin}}$-closed class.
	\end{definition}
	Aczel showed that we have the same $\HF$ if we replace $\Phi_\text{fin}$ to $\Phi_\text{f.e.}=\{(\{a\},a) \mid a\text{ is finitely enumerable}\}$ or $\Phi_\text{adj}=\{(\{a,b\},a\cup\{b\}) \mid a,b\in V\}$.
	Moreover, Aczel showed the following facts:
	\begin{proposition}\label{Proposition:InductionHF} \pushQED{\qed}
		\begin{enumerate}
			\item $\HF$ is transitive.
			\item If $\forall x\in \HF [(\forall y\in x \phi(y))\to \phi(x)]$, then $\forall x\in \HF \phi(x)$.
			\item If $\forall x,y\in \HF [(\forall u\in x\forall v\in y \phi(u,v))\to \phi(x,y)]$, then $\forall x,y\in \HF \phi(x,y)$. \qedhere 
		\end{enumerate}
	\end{proposition}
	
	\begin{proposition}\label{Proposition:DiscreteHF} \pushQED{\qed}
		$=$ and $\in$ is discrete over $\HF$, that is, $\forall x,y\in \HF (x=y\lor x\neq y)$ and $\forall x,y\in \HF (x\in y\lor x\notin y)$. \qedhere 
	\end{proposition}
	See Proposition 10.9 and 10.10 of \cite{Aczel2013} for the proof.

	We may ask whether $\HF$ is a set, as $\mathsf{ZF}$ proves $V_\omega$ is a set. The answer is affirmative, and we will prove it by constructing a hierarchy of $\HF$:
	\begin{definition}
		Let $\Dec(A)$ be a set of all decidable subsets of $A$:
		\begin{equation}
		\Dec(A) = \{B\subseteq A \mid \forall x\in A [ x\in B\lor x\notin B]\}.
		\end{equation}
		Define $D_n$ recursively as follows: $D_0=\varnothing$ and $D_{n+1}=\Dec(D_n)$.
	\end{definition}
	
	\begin{lemma}
		$\langle D_n\mid n\in\omega\rangle$ is a strictly increasing sequence of transitive discrete sets of $\HF$.
	\end{lemma}
	\begin{proof}
		It is easy to see that $D_n\subsetneq D_{n+1}$ by induction on $n$.
		For transitivity, observe that $x\in D_{n+1}$ implies $x\subseteq D_n\subseteq D_{n+1}$. Moreover, $\Dec(X)$ and ${{}^X}2$ have the same cardinality, so we can see each $D_n$ is finite.
		
		For discreteness of $D_n$, assume that $x,y\in D_{n+1}=\Dec(D_n)$. Then the formula $z\in x$ and $z\in y$ is decidable over $D_n$. Hence the formula $z\in x\lr z\in y$ is also decidable over $D_n$.
		Since $D_n$ is finite, we have 
		\begin{equation}
		\exists z\in D_n \lnot (z\in x\lr z\in y) \lor \lnot[\exists z\in D_n \lnot (z\in x\lr z\in y)],
		\end{equation}
		which implies $x\neq y \lor x=y$.
		
		It remains to show that $D_n\in \HF$. It follows from finiteness of $D_n$ and $D_n\subseteq\HF$ that will be shown by induction on $n$: if $D_n\subseteq\HF$ and $x\in D_{n+1}$, then $x\subseteq\HF$ is finite by \autoref{Lemma:DecidableFinite}. Hence $x\in \HF$ and we have $D_{n+1}\subseteq \HF$.
	\end{proof}
	
	\begin{theorem}
		$\HF = \bigcup_{n\in\omega} D_n$. Especially, $\HF$ is a set provided if $\omega$ is a set.
	\end{theorem}
	\begin{proof}
		Since $D_n\in \HF$ for all $n$ and $\HF$ is transitive, we have $\bigcup_{n\in\omega} D_n\subseteq\HF$.
		We will show that $\forall x\in\HF (x\in\bigcup_{n\in\omega}D_n)$ by applying \autoref{Proposition:InductionHF} to obtain the remaining inclusion.
		Assume that $x\subseteq \bigcup_{n\in\omega}D_n$ holds.
		Since $x$ is finite, we can find $n$ such that $x\subseteq D_n$. Choose a bijection $f\colon m\to x$ for some $m\in\omega$.
		We can see that for $y\in D_n$, $y\in x$ if and only if $\exists k<m (f(k)=y)$, which is a decidable formula since the quantifier is bounded by a natural number and $D_n$ is discrete. Hence $x\in \Dec(D_n)=D_{n+1}$.
	\end{proof}
	
	In classical world, $V_\omega$ is a model of finitary $\mathsf{ZF}$, which is bi-interpretable with $\mathsf{PA}$. Since $\mathsf{CZF^{fin}}$ is bi-interpretable with $\mathsf{HA}$ and the classical $V_\omega$ is the set of all hereditarily finite sets, we may ask  $\HF$ is a model of finitary $\mathsf{CZF}$. The following theorem shows the answer is affirmative:
	\begin{theorem}
		$\HF$ satisfies $\mathsf{CZF^{fin}}$, that is, if $\sigma$ is an axiom of $\mathsf{CZF^{fin}}$, then the relativization $\sigma^\HF$ of $\sigma$ holds.
	\end{theorem}
	\begin{proof}
		It suffices to show that $\sigma^\HF$ holds for all axioms of $\mathbb{T}$.
		Extensionality and Set Induction follow from \autoref{Proposition:InductionHF}. Moreover, $\HF$ is closed under the operation $x,y\mapsto x\cup \{y\}$, and this proves $\HF$ satisfies Pairing and Binary union. (For Binary union, we need to use the induction on the size of sets.) Since every $x\in\HF$ is finite, we can see $\bigcup x\in \HF$ by induction on the size of $x$. It shows $\HF$ satisfies Union.
		
		For Binary Intersection, let $x,y\in \HF$. Take $n\in\omega$ such that $x,y\in D_n$. By \autoref{Proposition:DiscreteHF}, both $z\in x$ and $z\in y$ is decidable over $\HF$. Hence the set $x\cap y = \{z\in D_n\mid z\in x\land z\in y\}$ is also decidable, so $x\cap y\in D_{n+1}$.
		
		It remains to show that $V=\Fin$ is valid in $\HF$. We know that for each $x\in \HF$ there is $n\in\omega$ and a bijection $f\colon n\to x$.
		We can see that $f$ is a finite set and $f\subseteq \HF$.
		Hence $f\in \HF$, and $f$ witnesses finiteness of $x$ in $\HF$.
	\end{proof}
	
	We may further ask $\HF$ satisfies finitary $\mathsf{IZF}$, which is identical with $\mathsf{ZF^{fin}}$. We will see that the answer is negative in general.
	\begin{theorem}\label{Theorem:HFnotIZFfin}
		Working over $V[\mathcal{K}_0]$, the model of Kleene realizability, $\HF$ does not satisfy the full law of excluded middle.
	\end{theorem}
	\begin{proof}
		It is known by \cite{McCartyPhD} and \cite{Rathjen2003Realizability} that $V[\mathcal{K}_0]$ satisfies \emph{Church's thesis} $\mathsf{CT_0}$. By Proposition 4.3.4 of \cite{DalenTroelstra}, $\mathsf{CT_0}$ implies the following instance of the negation of weak excluded middle $\mathsf{WLEM}$ holds:
		\begin{equation}\label{Formula:InstanceWLEM}
		\lnot\forall x\in\omega [\lnot\exists y\in\omega T(x,x,y) \lor \lnot\lnot\exists y\in\omega T(x,x,y)],
		\end{equation}
		where $T$ is Kleene's $T$-predicate, which has a $\mathcal{E}_1$-definition over $\omega$. Observe that the function $\mathfrak{p}$ defined in \autoref{Section:ordinalT} yields a definable bijection from $\omega$ to $\HF$, which is accessible inside $\HF$.
		Replacing all $x$ and $y$ of \eqref{Formula:InstanceWLEM} to $\mathfrak{p}^{-1}(x)$ and $\mathfrak{p}^{-1}(y)$ provides
		the formula of the form
		\begin{equation}
		\lnot\forall x\in\HF [\lnot\exists y\in\HF \phi(x,y) \lor \lnot\lnot\exists y\in\HF \phi(x,y)],	
		\end{equation}
		where $\phi(x,y)$ is a formula whose quantifiers inside this formula is bounded by $\HF$. Hence we may regard $\phi(x,y)$ as a relativization $\psi^\HF(x,y)$ of some formula $\psi(x,y)$.
		Therefore, $\HF$ satisfies an instance of the negation of $\mathsf{WLEM}$.
	\end{proof}
	
	%
	
	
	\section{Remarks and Questions}
	We will finish this article with a philosophical remark and some questions.
	We pointed out that $\mathbb{T}$ is bi-interpretable with $\mathsf{HA}$, and the bi-interpretation also captures bi-interpretability between subtheories of $\mathbb{T}$ and $\mathsf{HA}$. Moreover, $\mathbb{T}$ is $\mathsf{CZF^{fin}}$, and the set of all hereditarily finite sets $\HF$ is a model of $\mathbb{T}$.
	On the other hand, the finitary $\mathsf{IZF}$ is just $\mathsf{ZF^{fin}}$, which is not bi-interpretable with $\mathsf{HA}$. Moreover, \autoref{Theorem:HFnotIZFfin} shows $\HF$ may not be a model of finitary $\mathsf{IZF}$, even though the background universe satisfies $\mathsf{IZF}$: In \autoref{Theorem:HFnotIZFfin}, the background universe is $V[\mathcal{K}_0]$. If $V$ satisfies $\mathsf{IZF}$, then so does $V[\mathcal{K}_0]$.
	These two facts could bolster the viewpoint that $\mathsf{CZF}$ is a more natural constructive counterpart of $\mathsf{ZF}$ than $\mathsf{IZF}$.
	
	Classically, the negation of the axiom of Infinity proves $V=\Fin$. We do not know whether this is possible constructively. Instead, we postulate $V=\Fin$ as an axiom. It is natural to ask whether the negation of Infinity proves $V=\Fin$. Since $\mathbb{T}$ is an extension of Tharp's quasi-intuitionistic set theory \cite{Tharp1971} without the axiom of Infinity, and it contains the principle $\mathsf{Ord\mhyphen Im}$ that states every set is an image of an ordinal. Since $V=\Fin$ proves $\mathsf{Ord\mhyphen Im}$, we may also ask the question whether we can obtain an implication under $\mathsf{Ord\mhyphen Im}$:
	\begin{question}
		Does the negation of the axiom of Infinity prove $V=\Fin$? Can we prove it with $\mathsf{Ord\mhyphen Im}$?
	\end{question}
	
	We gave a bi-interpretation between $\IE_n$ and $\SIE_n$, which are subtheories of $\mathsf{HA}$ and $\mathbb{T}$ respectively. Unfortunately, we do not know the set-theoretic counterpart of constructive $\mathsf{I\Delta_0+Exp}$. Pettigrew \cite{Pettigrew2008} characterized the set-theoretic equivalent of $\mathsf{I\Delta_0+Exp}$, which is derived from Mayberry's set theory called $\mathsf{EA}$. Unfortunately, the auther does not know how to characterize the constructive counterpart of Pettigrew's set theory, hence the following question is still open:
	\begin{question}
		Can we identify a constructive set theory that is bi-interpretable with constructive $\mathsf{I\Delta_0+Exp}$?
	\end{question}

	The author only noticed after finishing the draft that McCarty and Shapiro had planned to gave an online talk on the Logic Supergroup (the detail will appear in \cite{ShapiroMcCarty}), and they independently achieved some of the author's result. The author was aware of their talk by chance on 21 September, only five days before their talk. The author contacted McCarty to inform the author's research, and he kindly responded with their slide, the main source the author can check the detail of their research.
	They worked with Heyting arithmetic with symbols for primitive recursion functions and another constructive set theory called $\mathsf{SST}$, which comprises Extensionality, the existence of adjunction $x+y:=x\cup\{y\}$, and adjunction induction
	\begin{equation}
		[\phi(0)\land \forall x\forall y (y\notin x \land \phi(x)\land \phi(y) \to \phi(x+y))]
		\to \forall x \phi(x)
	\end{equation}
	what they called Set Induction. They proved that the expanded $\mathsf{HA}$ is bi-interpretable, or \emph{definitionally equivalent} according to McCarty and Shapiro, with $\mathsf{SST}$ expanded by adding function symbols for primitive recursive functions. They also present a variant of Heyting arithmetic called $\mathsf{HA_{BIT}}$, and showed that it is bi-interpretable with $\mathsf{SST}$ and the extended $\mathsf{HA}$. Therefore, the extended $\mathsf{HA}$ and $\mathsf{SST}$ are bi-interpretable with each other.
	
	Moreover, we can show that $\mathsf{SST}$ is identical with $\mathsf{CZF^{fin}}$. The hardest part is deriving all axioms of $\mathsf{CZF^{fin}}$ with Collection instead of Strong Collection from $\mathsf{SST}$, and this is done by \cite{ShapiroMcCarty}.
	Hence we can see the author's bi-interpretability result for $\mathsf{HA}$ coincides with that of McCarty and Shapiro up to the difference of the language and choice of axioms of theories.
	
	The reviewer pointed out that Rathjen \cite{Rathjen2008Natural} provided a fragment of $\mathsf{CZF}$, which has the same proof-theoretic strength with $\mathsf{HA}$. Rathjen called this theory $\mathsf{CZF^-}$, and is obtained by discarding Set Induction from $\mathsf{CZF}$ and employing Strong Infinity instead of usual Infinity. Especially, $\mathsf{CZF^-}$ proves Mathematical Induction Axiom Scheme for $\Delta_0$-formulas.
	
	Rathjen proved in \cite{Rathjen2008Natural} that $\mathsf{CZF^-}$ is $\Pi^0_2$-conservative over $\mathsf{HA}$, and he used a mixture of type-theoretic interpretation of set theory and realizability interpretation of type theory over a saturated model of $\mathsf{PA}$.
	Rathjen also claimed that one can establish a similar synthetic translation into the theory $\mathsf{PA}^r_\Omega$, which is a conservative extension of $\mathsf{PA}$. However, Rathjen's translation of $\mathsf{CZF^-}$ to $\mathsf{PA}^r_\Omega$ is not an interpretation in our sense because his translation remolds quantifiers and logical connectives. Thus we have the following question:
	
	\begin{question}
	    Is there an interpretation from $\mathsf{CZF^-}$ to $\mathsf{HA}$ (or equivalently, $\mathsf{CZF_{fin}}$)? 
	\end{question}

	\section*{Acknowledgement}
    This paper is the author's master's thesis.
	The author would like to thank the author's advisor Otto van Koert for improving the grammatical and style of this paper. The author is also grateful to the referee for pointing out grammatical errors and making helpful suggestions.
	The author also wants to express the gratitude to David Charles McCarty, Stewart Shapiro and Andr\'es E. Caicedo for their encouragement on the author's work.

	\overfullrule=0pt
	\printbibliography
	\nocite{*}
	
	
\end{document}